\documentclass[journal,twoside]{IEEEtran}
\usepackage{fancyvrb}
\usepackage{listings}
\usepackage[pdf]{pstricks}
\usepackage{cite}
\ifCLASSINFOpdf
  \usepackage[pdftex]{graphicx}
\else
  \usepackage[dvips]{graphicx}
\fi

\usepackage{amsmath,mathtools,amsfonts,amssymb,amscd,bm,cases}
\usepackage{enumitem}
\usepackage{hyperref}
\usepackage{amsmath,amssymb}
\usepackage{cancel}
\usepackage{graphics,psfrag}
\usepackage{empheq}
\usepackage{colortbl}

\usepackage{stmaryrd}
\makeatletter
\newcommand{\xMapsto}[2][]{\ext@arrow 0599{\Mapstofill@}{#1}{#2}}
\def\Mapstofill@{\arrowfill@{\Mapstochar\Relbar}\Relbar\Rightarrow}
\makeatother

\begin{document}
%~ \title{Covariant mappings optimize the parametric Discontinuous Galerkin time domain solution of Maxwell equations and shed light on complex stretching based Perfectly Matched Layers}
\title{A matrix-free Discontinuous Galerkin method\\ for the time dependent Maxwell equations\\ in unbounded domains}
\author{Bernard~Kapidani and Joachim Sch{\"o}berl
%~ \thanks{Manuscript received July 21, 2017.}% <-this % stops a space
\thanks{The authors are with the Institute for Analysis and Scientific Computing, Technische Universit{\"a}t Wien, Wiedner Hauptstr. 108--110, A-1040 Vienna, Austria (e-mail: bernard.kapidani@tuwien.ac.at).}% <-this % stops a space
}% <-this % stops a space
% The paper headers
% \markboth{IEEE Transactions on Antennas and Propagation,~Vol.~11, No.~4, December~2016}%
% {Codecasa \MakeLowercase{\textit{et al.}}: A Consistent and Conditionally Stable Time-Domain Explicit Scheme on Tetrahedral Meshes for Wave Propagation Problems in Lossy Materials}
%\markboth{IEEE TRANSACTIONS ON ANTENNAS AND PROPAGATION,~Vol.~XX, No.~X, Month~XXXX}{IEEE TRANSACTIONS ON ANTENNAS AND PROPAGATION,~Vol.~XX, No.~X, Month~XXXX}%
% make the title area
\maketitle

\vspace{-2cm}
% As a general rule, do not put math, special symbols or citations
% in the abstract or keywords.
\begin{abstract}
A Discontinuous Galerkin (DG) Finite Element Method (FEM) approach for the 3D time dependent Maxwell equations in unbounded domains is presented. The method (implemented in the FEM library NGsolve) is based on the covariant transformation of a modal orthogonal polynomial basis, originally defined on a reference simplex. The approach leads to an explicit time stepping scheme for which the mass matrix to be inverted is at most $d\!\times\!d$ block-diagonal (in $d\!=\!2,3$ spatial dimensions) while the matrix which discretizes the curl operators on the right-hand side of the system is a small reference matrix, independent from geometric properties of mesh elements. Furthermore, we show that the introduced optimizations are preserved when unbounded domains are also included in the formulation through a complex-stretching based approach.
\end{abstract}
% Note that keywords are not normally used for peerreview papers.
\IEEEpeerreviewmaketitle
\section{Introduction}
Technological advancements in the field of wideband devices and telecommunication techniques operating at frequencies in the MHz range and above have spurned a continued research effort, in the field of computational science, in improving numerical methods for the solution of the time dependent Maxwell equations (ME henceforth), i.e.
\begin{align}
 &\partial_t \bm{d}\left(\bm{r},t\right) = \nabla\!\times\! \bm{h}\left(\bm{r},t\right) - \mathbf{j}\left(\bm{r},t\right),   \label{eq:ampere_cont} \\
 &\partial_t \bm{b}\left(\bm{r},t\right) = - \nabla\!\times\! \bm{e}\left(\bm{r},t\right), \label{eq:faraday_cont}
\end{align}
\noindent supplemented with the usual (time-invariant) constitutive equations 
\begin{align*}
%\label{eq:eps_cont}
 &\bm{d}\left(\bm{r},t\right) = \varepsilon \left(\bm{r}\right) \bm{e}\left(\bm{r},t\right), \;\; \bm{b}\left(\bm{r},t\right) = \mu \left(\bm{r}\right) \bm{h}\left(\bm{r},t\right),
\end{align*}
\noindent and to be solved for a (possibly unbounded) domain $\Omega \subseteq \mathbb{R}^d$ and a finite time interval $[0,\mathrm{T}]$ with suitable initial conditions (we hereinafter neglect conduction currents for ease of presentation, only to reintroduce them at the appropriate point in the article). If an accurate numerical time domain solution of (\ref{eq:ampere_cont})-(\ref{eq:faraday_cont}) is available, valuable information on the simulated device or electromagnetic structure can be obtained for a many frequencies in one stroke by post-processing the results with Fast Fourier Transform (FFT) computations.

In this endeavour, a fair amount of literature has been produced in finite differences based approximation, starting from the seminal paper of Yee \cite{yee_numerical_1966}. The monograph on the subject by Taflove \cite{a_taflove_computational_2000} is a comprehensive text on all related subsequent advancements. With some delay, it was shown\cite{bossavit_solving_1990,jin-fa_lee_whitney_1995} that the Finite Element Method (FEM) was also a viable path for the same electromagnetic wave propagation problem, namely, by using basis functions originally introduced by Nedelec in \cite{nedelec_mixed_1980} (which were actually a re-discovery of functions introduced by Whitney in 1958 \cite{whitney_geometric_2015}). These so-called \emph{edge elements} are piecewise-polynomial vector-valued basis functions on simplicial meshes, having continuous tangential components and discontinuous normal ones on boundaries between neighbouring tetrahedra in the spatial discretization. This partial continuity properties allow a faithful discretization of (\ref{eq:ampere_cont})-(\ref{eq:faraday_cont}) in the general setting of problems with discontinuous material coefficients. In later years, extensions of the FEM approach have been introduced with higher polynomial order in the approximating functions\cite{graglia_higher_1997} and with hierarchical structure for the bases\cite{webb_hierarchal_1999,schoberl_high_2005} (which fact allows to locally control the polynomial approximation degree on each element in the mesh).

Unfortunately, the flexibility in modelling complex geometries given by meshing with tetrahedra and triangles, and the arbitrarily high polynomial degree in the approximation, are offset by the structure of the resulting algebraic system to be solved. It follows from the support of edge basis functions that the mass matrix obtained by Galerkin testing of the Ampere--Maxwell and Faraday equations, is a sparse and banded one, but not block-diagonal. This fact has impactful implications when the time dependent problem needs to be solved, since the mass matrix needs to be inverted by means of an iterative solver at every time step. Any direct factorization approach to compute such inverse is hardly scalable, since in problems of practical interest the number of unknowns is easily in the order of millions.

To overcome this issue, the use of fully non-conforming finite elements has been advocated, extending the Discontinuous Galerkin (DG) approach to the time-dependent Maxwell problem (see \cite{hesthaven_nodal_2002}, \cite{jin_finite_2014}). The basic idea in this setting is to approximate the unknown electric and magnetic fields in each finite element with a local basis, complete up to some polynomial degree $p$, without regard for inter-element continuity constraints on basis functions. The resulting algebraic system then presents naturally block-diagonal mass matrices (one block per finite element).

This approach has also undesired side-effects. Firstly, additional care must be taken on the right-hand side (r.h.s.) of the equations, as integration by part on each element will lead to the appearance of inter-element boundary integrals in which multi-valued fields are involved. These jumps in tangential components of fields must be penalized, which fact leads usually to numerical schemes which artificially dissipate electromagnetic energy.
Furthermore, in practice, very high local polynomial order in the approximating functions needs to be employed, where the global set of unknowns will be equal to the Cartesian product of local element unknowns, i.e. there is no condensation of degrees of freedom (DoFs) when the algebraic system on the global mesh is assembled, as is instead the case with edge elements.
It is therefore paramount to make implementation choices which reduce the time and memory costs of the DG formulation. 

In the following, we describe an energy-conserving DG method for the time-dependent ME in which the vector-valued basis functions are based on orthogonal polynomials on a reference simplex (defined in Section \ref{sec:spaces}). Furthermore, appropriate application of the covariant and contravariant transformation rules of vector fields from the reference to physical finite elements result in a discrete approximation of the curl operator which does not dependent on the geometry of the elements and can be therefore computed and stored in memory as a relatively small matrix for the reference element (see Section \ref{sec:weak_form}). Most importantly, the scheme will be applied for the first time to the solution of problems in unbounded domains: for this end we employ complex coordinate-stretching based Perfectly Matched Layers (PML)\cite{teixeira_complex_2000} which, when used in the time domain, lead to additional vector unknowns to be solved for in the formulation. It will be shown, in Section \ref{sec:pml_tf}, that the transformation rules underlying the PML are in the same family of the already exploited mappings between reference and physical finite elements, and thus preserve (and further exploit) the introduced computational optimizations. Numerical experiments validate the proposed method in Section \ref{sec:num_res} and conclusions are drawn in Section \ref{sec:cr}.

\section{Discontinuous approximation spaces}\label{sec:spaces}
We will start from the definition of basis functions for the approximation of the unknown fields on a reference element ${\hat{T}_3}$, as common practice in the standard FEM. In three spatial dimensions, our reference domain will be the tetrahedron defined as the set ${\hat{T}_d} = \{ \hat{\bm{r}}(\hat{x},\hat{y},\hat{z})$ such that $\hat{x},\hat{y},\hat{z} \geq 0$ and $\hat{x}+\hat{y}+\hat{z} \leq 1$. In two dimensions, we accordingly use the triangle ${\hat{T}_2} = \{ \hat{\bm{r}}=(\hat{x},\hat{y}) \; : \; \hat{x},\hat{y} \geq 0 \,\wedge\, \hat{x}+\hat{y} \leq 1\}$.

We choose a set of mutually orthogonal basis functions with support ${\hat{T}_d}$, where the orthogonality condition is here meant in the space of square integrable functions $L^2({\hat{T}_d})$. For space dimension $d$=3, the following modal basis is available:
\begin{align}
  \begin{split}
    \hat{\varphi}_{ijk} (\hat{\boldsymbol{r}}) =
    &P_i^{(0,0)} \left(\frac{2\hat{z}}{1-\hat{x}-\hat{y}} - 1\right) \left(1-\hat{x}-\hat{y}\right)^i\times\\
    &P_j^{(2i+1,0)} \left(\frac{2\hat{y}}{1-\hat{x}} - 1\right)\times\\
    &P_k^{(2i+1,0)} \left(2\hat{x} - 1\right),
  \end{split}
  \label{eq:dubiner_functions}
\end{align}
\noindent with $i,j,k$ non-negative integers such that $i+j+k \leq p$ and where $P_n^{(\alpha,\beta)} (w)$ are the classical Jacobi polynomials, with the generating formula:
\begin{align*}
\begin{split}
  P_n^{(\alpha,\beta)} (w) = &\frac{(-1)^n}{2^n n!} (1 - w)^{-\alpha} (1 + w)^{-\beta} \times\\
                          &\frac{\mathrm{d}^n}{\mathrm{d}w^n}\biggl( (1-w)^\alpha (1+w)^\beta (1-w)^n (1+w)^n\biggr),
\end{split}
\end{align*}
\noindent which we recall are $\mathrm{L}^2$--orthogonal with respect to the weight $(1-x)^\alpha (1+x)^\beta$ on the interval $[-1,1]$. Furthermore the special case $\alpha = \beta = 0$ yields the well-known Legendre polynomials.
%% \end{itemize}

The set of polynomials defined in (\ref{eq:dubiner_functions}), for which an equivalent in the case $d$=2 can be found in \cite{dubiner_spectral_1991}, satisfies then
\begin{align}
  \left( \hat{\varphi}_{ijk}(\hat{\bm{r}}),\hat{\varphi}_{i'j'k'}(\hat{\bm{r}})\right)_T &:= 
  \int_{{\hat{T}_d}} \hat{\varphi}_{ijk}(\hat{\bm{r}}) \hat{\varphi}_{i'j'k'}(\hat{\bm{r}}) \mathrm{d}\hat{\bm{r}} = \nonumber\\
  & = \delta_{ii'} \delta_{jj'}, \delta_{kk'},\label{eq:ortho_cond}
\end{align}
\noindent where $\delta$ is the Kronecker tensor.
We remark that the orthogonality of polynomials in (\ref{eq:dubiner_functions}) is achieved by construction from the orthogonality of a tensor product of (shifted and scaled) univariate Jacobi polynomials (in $\hat{x}$, $\hat{y}$, $\hat{z}$ respectively) on the reference cube $\hat{Q}\!=\![0,1]\!\times\![0,1]\!\times\![0,1]$: the domain ${\hat{T}_d}$ can be in fact obtained from $\hat{Q}$ via a Duffy transformation\cite{duffy_quadrature_1982}, which (for $d$=3) consists in collapsing the upper facet of $\hat{Q}$ into the point with coordinates $(0,0,1)^{\mathrm{T}}$. 
Analogous considerations apply for the case $d=2$ (with the appropriate Duffy transformation from a square to a triangle).

We want to approximate vector fields on ${\hat{T}_d}$ within $\mathbb{P}^p({\hat{T}_d};\mathbb{R}^d)$, i.e. the space of polynomials of degree lower or equal to $p$ on ${\hat{T}_d}$.
By taking 
\begin{align*}N_d^p = {p+d\choose d}\end{align*} 
\noindent distinct orthogonal $\varphi_{ijk}(\hat{\boldsymbol{r}})$ for each component, we construct a fully vectorial basis:
\begin{align*}
  & {\hat{\boldsymbol{\varphi}}}_{j\ell}(\hat{\bm{r}}) := {\hat{\varphi}}_{\ell}(\hat{\bm{r}}) \hat{\mathbf{e}}_j,
\end{align*}
%% \begin{align*}
%%   & \hat{\boldsymbol{e}}\left(\hat{\boldsymbol{r}},t\right)  = \sum\limits_{\ell=1}^{\ell=N} \sum\limits_{j=1}^{j=3} u_{j\ell}(t) \varphi_{\ell}(\hat{\mathbf{r}}) \hat{\mathbf{e}}_j \;\; \forall \hat{\boldsymbol{r}} \in {\hat{T}_d}, \\
%%   & \hat{\boldsymbol{h}}\left(\hat{\boldsymbol{r}},t\right)  = \sum\limits_{\ell=1}^{\ell=N} \sum\limits_{j=1}^{j=3} f_{j\ell}(t) \varphi_{\ell}(\hat{\mathbf{r}}) \hat{\mathbf{e}}_j \;\; \forall \hat{\boldsymbol{r}} \in {\hat{T}_d},
%% \end{align*}
\noindent where $\hat{\mathbf{e}}_j$ is the $j$-th Cartesian coordinate unit vector and $\ell\in\{1,2,...,N_d^p\}$ can be substituted for the three indices $i$,$j$,$k$ without ambiguity. With the integrand of (\ref{eq:ortho_cond}) turned into a dot product, one gets a fully orthogonal vector-valued basis:
\begin{align*}
  \left(  \hat{\boldsymbol{\varphi}}_{j\ell}(\hat{\bm{r}}) , \hat{\boldsymbol{\varphi}}_{j'\ell'}(\hat{\bm{r}}) \right)_{{\hat{T}_d}} := \int_{{\hat{T}_d}} \hat{\boldsymbol{\varphi}}_{j\ell}(\hat{\bm{r}}) \cdot \hat{\boldsymbol{\varphi}}_{j'\ell'}(\hat{\bm{r}}) \mathrm{d}\hat{\bm{r}} = \delta_{jj'} \delta_{\ell\ell'}.%\label{eq:ortho_cond_vector}
\end{align*}

We assume now a conforming partition of our spatial domain of interest $\Omega$ into simplexes, which we denote with $\mathcal{T}_\Omega$. Any simplex $T\subset\mathcal{T}_\Omega$ can be generated from ${\hat{T}_d}$ via an affine invertible mapping ${\Phi}_{T}: {\hat{T}_d} \mapsto T,$ such that
\begin{align}
  & \bm{r} = {\Phi}_{T}(\hat{\bm{r}}) = \mathbf{A}_T \hat{\bm{r}} + \mathbf{b}_T,\label{eq:mapping_template}
\end{align}
\noindent where $\mathbf{A}_T\in\mathbb{R}^{d\times d}$ and $\mathbf{b}_T\in\mathbb{R}^d$. Furthermore, we denote with $\hat{\bm{r}}$ position vectors in the local coordinates of ${\hat{T}_d}$ and with $\bm{r}$ points in the global ones, with $\bm{r} \in T \subset \mathcal{T}_\Omega$.

The geometric mapping in (\ref{eq:mapping_template}) induces the transformation rules for scalar and vector fields expanded in the local and global basis coordinates. For example, the one-to-one mapping from scalar-valued functions $\hat{f}(\hat{\boldsymbol{r}}) \in \mathbb{P}^{p+1}({\hat{T}_d};\mathbb{R})$ to functions $f(\boldsymbol{r})\in\mathbb{P}^{p+1}(T;\mathbb{R})$ is given by:
\begin{align*}
f(\boldsymbol{r}) = f \circ {\Phi}_{T} (\hat{\boldsymbol{r}}) = \hat{f} (\hat{\boldsymbol{r}}) \implies 
  \hat{f} (\hat{\boldsymbol{r}}) = \hat{f} \circ {\Phi}_{T}^{-1} (\boldsymbol{r}).%\label{eq:composition}
\end{align*}

If we also need to compute derivatives of $f(\bm{r})$, it is easy to show, by virtue of the chain-rule, that the gradient of $f(\bm{r})$ in the global coordinates $\nabla f (\bm{r}) \in (\mathbb{P}^p(T);\mathbb{R}^d)$ is given (again for $d=3$) by
\begin{align} 
  & \nabla f (\boldsymbol{r}) =
  \begin{pmatrix*}[c]
    \vspace{0mm}\partial_{\hat{x}}{x}  & \partial_{\hat{y}}{x} & \partial_{\hat{z}}{x} \\
    \vspace{0mm}\partial_{\hat{x}}{y}  & \partial_{\hat{y}}{y} & \partial_{\hat{z}}{y}\\
    \vspace{0mm}\partial_{\hat{x}}{z}  & \partial_{\hat{y}}{z} & \partial_{\hat{z}}{z}
  \end{pmatrix*}\!
\begin{pmatrix*}
 \vspace{0mm} \partial_{\hat{x}} \hat{f} \\
  \vspace{0mm} \partial_{\hat{y}} \hat{f} \\
  \vspace{0mm} \partial_{\hat{z}} \hat{f} \end{pmatrix*} = 
  \mathbf{A}_T^{-\mathrm{T}} \hat{\nabla} \hat{f}(\hat{\bm{r}}),\label{eq:cov_transf}
\end{align}
\noindent where $(\cdot)^{-\mathrm{T}}$ denotes the inverse-transpose, $\mathbf{A}_T$ is the Jacobian matrix of the mapping ${\Phi}_{T}$, and symbol $\hat{\nabla}$ denotes the gradient computed with respect to the local $\hat{x}$, $\hat{y}$, and $\hat{z}$ coordinates on ${\hat{T}_d}$.

The transformation matrix $\mathbf{A}_T^{-\mathrm{T}}$ in Eq. (\ref{eq:cov_transf}) represents a covariant transformation. When using conforming Finite Element discretizations for the ME, it is well known that the covariant transformation is the appropriate one for mapping basis functions defined on a reference element to their counterparts on the physical element (see for example \cite{monk_finite_2003, matuszyk_parametric_2013}). 
The reason for this lies in the fact that $\mathbf{A}_T^{-\mathrm{T}}$ transforms locally curl-free fields on ${\hat{T}_d}$ into curl-free fields on $T$ (recall $\nabla\!\times\!\nabla f = 0$ $\forall f$, also known in literature as De Rham sequence properties\cite{monk_finite_2003}).
To make this more clear one can reason in the framework of a standard Helmholtz decomposition for any vector field $\hat{\bm{v}}\in\mathbb{P}^{p}({\hat{T}_d};\mathbb{R}^d)$:
\begin{align*}
& \hat{\bm{v}} = \hat{\nabla}\hat{f} + \hat{\bm{g}},
\end{align*}
where also $\hat{\bm{g}}\in\mathbb{P}^{p}({\hat{T}_d};\mathbb{R}^2)$ and 
\begin{align*}\left(\hat{\bm{g}},\hat{\nabla}\hat{f}\right)_{{\hat{T}_d}} = 0 \implies \left(\bm{g},\nabla f\right)_{T} = 0,\end{align*}
\noindent where $\bm{g} = \mathbf{A}_T^{-\mathrm{T}} \hat{\bm{g}}$.
Even if we are dealing with completely discontinuous basis functions, we use the rule in (\ref{eq:cov_transf}) and, starting from any member of the corresponding basis on the reference simplex, we define, in the global coordinates, functions
\begin{align}
\boldsymbol{\varphi}_{j\ell}^T (\bm{r}) = 
\begin{cases}\mathbf{A}_T^{-\mathrm{T}} \hat{\boldsymbol{\varphi}}_{j\ell} ({\Phi}_{T}^{-1}({\bm{r}})), & \forall \bm{r}\in T \\
0 & \text{otherwise}, \end{cases}\label{eq:mapped_eh}
\end{align}
\noindent where it is important to note that, after mapping, in general
\begin{align*}
& \left(\boldsymbol{\varphi}_{j\ell}^T (\bm{r}),\boldsymbol{\varphi}_{j'\ell}^T (\bm{r})\right) = (\mathbb{A}^{-1}\mathbb{A}^{-T})_{jj'} \delta_{\ell\ell'},
\end{align*}
\noindent i.e. we lose some degree of orthogonality in the functional set, due to the physical element being in general a tetrahedron with no right-angle corners. Nevertheless the resulting inner product (a.k.a. mass) matrices are at worst $d\!\times\!d$ block-diagonal. This seems at first glance a drawback in terms of performance, yet its appeal will be made clear in Section \ref{sec:weak_form}.

Once we have established (\ref{eq:mapped_eh}) as our basis functions for $\bm{e}(\bm{r},t)$ and $\bm{h}(\bm{r},t)$ on each $T\in\mathcal{T}_\Omega$, we remark we will need to compute the curl of the local approximation of the fields. The chain rule and some tensor algebra lead us to derive the following transformation rule:
\begin{align}
\nabla\!\times\! \boldsymbol{\varphi}_{j\ell}^T (\bm{r}) = |\mathbf{A}_T|^{-1} \mathbf{A}_T\,\hat{\nabla} \times \hat{\boldsymbol{\varphi}}_{j\ell} (\,{\Phi}_{T}^{-1}(\bm{r})\,),\label{eq:piola_transf}
\end{align}
\noindent where $|\mathbf{A}_T|$ is the determinant of the Jacobian matrix $\mathbf{A}_T$. The transformation in (\ref{eq:piola_transf}) is known as the contravariant or \emph{Piola} mapping (in the context of continuum mechanics, where it is used to enforce mass conservation). Similarly to the covariant transformation, it maps the kernel of the divergence operator on the reference element to the the kernel of the divergence operator on the physical element (recall $\nabla\cdot\nabla\!\times\!\bm{v} = 0$ $\forall \bm{v}$). If (\ref{eq:cov_transf}) is used to define basis functions on the physical element, (\ref{eq:piola_transf}) will be needed to compute the local curl on each element in the discrete system.

\section{Discontinuous Galerkin formulation}\label{sec:weak_form}
In bottom-up fashion, once we have established the local approximation spaces, we can proceed to the DG solution for the ME. They key operational difference when using fully discontinuous basis functions with respect to the classical conforming FEM approach is that in the former one has to start by introducing the mesh (our aforementioned triangulation $\mathcal{T}_\Omega$ of the domain $\Omega$), and only subsequently one can proceed with Galerkin testing of the partial differential equations, which is carried out integrating locally (element-by-element) on each $T$.
We require (\ref{eq:ampere_cont})--(\ref{eq:faraday_cont}) to hold in the weak sense, i.e.
\begin{align}
&\sum\limits_{T \in \mathcal{T}_\Omega} \left(\,\varepsilon \partial_t\bm{e} , \boldsymbol{\varphi}_{j\ell}^T\right)_T =
   \sum\limits_{T \in \mathcal{T}_\Omega} \left( \nabla\!\times\! \bm{h}, \boldsymbol{\varphi}_{j\ell}^T\right)_T, \label{eq:weak1}\\
&\sum\limits_{T \in \mathcal{T}_\Omega} \left(\, \mu \partial_t\bm{h} , \boldsymbol{\varphi}_{j\ell}^T\right)_T =
 \sum\limits_{T \in \mathcal{T}_\Omega} \left( -\nabla\!\times\! \bm{e}, \boldsymbol{\varphi}_{j\ell}^T\right)_T, \label{eq:weak2}
\end{align}
\noindent must hold for all $\boldsymbol{\varphi}_{j\ell}^T$, $\ell \in 1,2,\dots, N_3^p$ and $j\in\{1,2,3\}$ on each $T\in\mathcal{T}_\Omega$, where we will omit $\bm{r}$ and $t$ dependence for readability where no confusion can arise.
A solution to the weak problem is accordingly sought in the same finite dimensional trial space, i.e. with the expansions:
\begin{align}
& \boldsymbol{e}\left(\boldsymbol{r},t\right)\!\mid_T  = \sum\limits_{\ell=1}^{\ell=N_d^{p+1}}\sum\limits_{j=1}^{j=3} u_{j\ell}(t) \boldsymbol{\varphi}_{j\ell}^T(\boldsymbol{r}), \label{eq:triale}\\
& \boldsymbol{h}\left(\boldsymbol{r},t\right)\!\mid_T = \sum\limits_{\ell=1}^{\ell=N_d^p} \sum\limits_{j=1}^{j=3} f_{j\ell}(t) \boldsymbol{\varphi}_{j\ell}^T(\boldsymbol{r}),\label{eq:trialh}
\end{align}
\noindent making (\ref{eq:weak1})--(\ref{eq:weak2}) a square system of ordinary differential equations. Indeed, since the time variable has not been discretized yet, degrees of freedom $u_{j\ell}(t)$ and $f_{j\ell}(t)$ explicitly encode the time dependence.

Take now any pair of trial and test functions: for the mass matrix entry associated to their inner product on the left-hand side (l.h.s.) of Eq. (\ref{eq:weak1})--(\ref{eq:weak2}), it holds
\begin{align*}
&\left( \varepsilon \boldsymbol{\varphi}_{ j\ell}^T, \boldsymbol{\varphi}_{km}^T\right)_T =
    \left( |\mathbf{A}_T| \varepsilon \mathbf{A}_T^{-\mathrm{T}} \hat{\boldsymbol{\varphi}}_{ j\ell} ,  \mathbf{A}_T^{-T}\hat{\boldsymbol{\varphi}}_{km}\right)_{{\hat{T}_d}},
\end{align*}
\noindent with the determinant factor arising from the change of variable in the integration. Furthermore, again for any pair of trial and test functions and for the bilinear forms on the r.h.s. of (\ref{eq:weak1})--(\ref{eq:weak2}), it holds
\begin{align}
 \left( \nabla\!\times\! \boldsymbol{\varphi}_{j\ell}^T, \boldsymbol{\varphi}_{km}^T\right)_T &=
 \left(|\mathbf{A}_T|^{-1} \mathbf{A}_T \hat{\nabla} \times \hat{\boldsymbol{\varphi}}_{ j\ell}, |\mathbf{A}_T|\mathbf{A}_T^{-\mathrm{T}}\hat{\boldsymbol{\varphi}}_{km} \right)_{{\hat{T}_d}}\!= \nonumber \\
&= \left(\hat{\nabla}\!\times\!\hat{\boldsymbol{\varphi}}_{ j\ell} , \hat{\boldsymbol{\varphi}}_{km} \right)_{{\hat{T}_d}},\label{eq:matrix_free_curl}
\end{align}
\noindent where (\ref{eq:piola_transf}) has been used to map curls of covariant vector fields. Eq. (\ref{eq:matrix_free_curl}) shows that the discrete bilinear forms on the r.h.s. of the equations, which discretize the local curl operator, are independent of mesh geometry, and can thus be computed once and for all (at the start of computation) as a small local matrix on the reference element and re-used for all DoFs in the mesh. Even for high polynomial orders in the trial and test functions, the curl operator matrix will thus grow moderately in size and, by virtue of exploiting cache locality, matrix-vector multiplications involved in the application of curl operators will be considerably sped-up.

Unfortunately the formulation is not usable as it stands, since no coupling between neighbouring elements occurs yet. To fix the issue we first proceed by formally integrating by parts the r.h.s. of the first weak equation
\begin{align}
&\sum\limits_{T \in \mathcal{T}_\Omega} \left( \varepsilon\partial_t\bm{e} , \boldsymbol{\varphi}_{j\ell}^T\right)_T =
\sum\limits_{T \in \mathcal{T}_\Omega}  \left( \nabla\!\times\! \bm{h} , \boldsymbol{\varphi}_{j\ell}^T\right)_T = \nonumber \\
&=\sum\limits_{T \in \mathcal{T}_\Omega}\hspace{-0mm}\left(\left( \boldsymbol{h} , \nabla\!\times\! \boldsymbol{\varphi}_{j\ell}^T\right)_T\!-\!\sum\limits_{F \in \partial{T}} \left( \{\!\{ \bm{h}\}\!\}_{\hat{\bm{t}}_F} , \boldsymbol{\varphi}_{j\ell}^T\right)_{F}\right), \label{eq:uwc}
\end{align}
\noindent where $\hat{\bm{t}}_F$ is the unit tangent vector on facet $F$ and $\{\!\{\cdot\}\!\}_{\hat{\bm{t}}_F}$ denotes the tangential averaging operator
\begin{align*}
&\{\!\{ \bm{h}\}\!\}_{\hat{\bm{t}}_F}  := \left(\frac{\bm{h}|_T+\bm{h}|_{\Omega\setminus T}}{2}\right)\bigg|_F\hspace{-3mm}\times\hat{\bm{n}}_F,\;\;\forall F\in\partial T, T\in\mathcal{T}_\Omega,
\end{align*}
\noindent with $\bm{h}|_{\Omega\setminus T}$ denoting the field value approximation taken from the neighbouring elements to $T$, i.e. each element sharing a mesh facet $F \subset \partial{T}$. Furthermore, $\hat{\boldsymbol{n}}_F$ is the (always outwards pointing) normal vector to the element boundary, itself denoted by $\partial{T}$. These so-called central (i.e. arithmetically averaged) numerical fluxes allow to preserve electromagnetic energy conservation in the ME, if we also add the following consistent term to the second weak equation:
\begin{align}
\sum\limits_{T \in \mathcal{T}_\Omega} \left( \mu \partial_t\boldsymbol{h} , {\boldsymbol{\varphi}}_{j\ell}^T\right)_T &= 
\sum\limits_{T \in \mathcal{T}_\Omega}\!\sum\limits_{F \in \partial{T}}\!\left( \{\!\{ \bm{e}\}\!\}_{\hat{\bm{t}}_F} \!-\!\boldsymbol{e} \times \hat{\boldsymbol{n}} , {\boldsymbol{\varphi}}_{j\ell}^T\right)_{F}+\nonumber\\
&-\sum\limits_{T \in \mathcal{T}_\Omega}\left( \nabla \!\times \!\boldsymbol{e} , {\boldsymbol{\varphi}}_{j\ell}^T\right)_T, \label{eq:cwc}
\end{align}
\noindent where all differential operators act on the field which is approximated with the higher polynomial degree in both weak equations.
%It is a matter of straightforward checking that, once the unknowns are expanded in the bases of (\ref{eq:triale})--(\ref{eq:trialh}), the resulting semi-discrete system retains the following hyperbolic form
%\begin{align*}
%             \begin{pmatrix*}[c]
%              \mathbf{M}_{\varepsilon} &  \mathbf{0} \\
%              \mathbf{0} &  \mathbf{M}_\mu
%             \end{pmatrix*}\frac{\mathrm{d}}{\mathrm{d}t}
%             \begin{pmatrix*}[c]
%              {\mathbf{e}} \\
%              {\mathbf{h}} 
%             \end{pmatrix*} =
%             \begin{pmatrix*}[c]
%              \mathbf{0} & \mathbf{C} \\
%              -\mathbf{C}_{}^{\mathrm{T}} & \mathbf{0} &
%             \end{pmatrix*}
%             \begin{pmatrix*}[c]
%              {\mathbf{e}} \\
%              {\mathbf{h}}
%             \end{pmatrix*},
%\end{align*}
%\noindent where $\mathbf{e}$, $\mathbf{h}$ are the column vectors containing the DoFs $u_{j\ell}(t)$ and $f_{j\ell}(t)$, respectively. 
To make the scheme also practically free of spurious solutions, we employ further the symmetric interior penalty (SIP) approach\cite{perugia_stabilized_2002}, which in the time-domain amounts to treating the tangential component jump of one field as a magnetic surface current, to which a contribution to the magnetic field corresponds, whose domain of definition is the skeleton of $\mathcal{T}_\Omega$:
\begin{align}
& \mu\partial_t \bm{h}_s = \alpha [[\bm{e}]]_{\hat{\bm{t}}_F}\;\;\forall F\in\partial T, \forall T\in\mathcal{T}_\Omega\label{eq:jmp_var_def}
\end{align}
\noindent where $\alpha>0$ is a real-valued penalty parameter depending on mesh size and $p$ and 
\begin{align*}
&[[\bm{e}]]_{\hat{\bm{t}}_F} := \{\!\{ \bm{e}\}\!\}_{\hat{\bm{t}}_F}\!-\!\bm{e}\bigg|_{\Omega\setminus T}\hspace{-6mm}\times\!\hat{\bm{n}}_F,\;\;\forall F\in\partial T, T\in\mathcal{T}_\Omega,
\end{align*}
\noindent is the tangential jump operator. Sufficient lower bounds for the value of $\alpha$ can be found in the literature (see \cite{warburton_role_2006,geevers_sharp_2017}).
The auxiliary unknown $\bm{h}_s(\bm{r},t)$ enters the formulation as a jump penalization term in the first equation (\ref{eq:uwc}), which becomes
\begin{align}
\sum\limits_{T\in\mathcal{T}_\Omega}\left(\partial_t\bm{e},\boldsymbol{\varphi}_{j\ell}^T\right)_T&=
\sum\limits_{T \in \mathcal{T}_\Omega}\left( \boldsymbol{h} , \nabla\!\times\! \boldsymbol{\varphi}_{j\ell}^T\,\right)_T+\nonumber\\
&-\sum\limits_{T \in \mathcal{T}_\Omega}\sum\limits_{F \in \partial{T}} \left( \{\!\{ \bm{h}\}\!\}_{\hat{\bm{t}}_F}  , \boldsymbol{\varphi}_{j\ell}^T\right)_{F} + \nonumber\\
&- \sum\limits_{T \in \mathcal{T}_\Omega} \sum\limits_{F \in \partial{T}} \left(\bm{h}_s\!\times\!\hat{\bm{n}}_F, \boldsymbol{\varphi}_{j\ell}^T\right)_{F},\label{eq:1st_weak_equation}
\end{align}
\noindent for any $\boldsymbol{\varphi}_{j\ell}^T$, and where all fields appearing in boundary integrals are single-valued and the added penalization term is consistent, since the tangential jumps vanish, $\forall F$ in the skeleton of the mesh.
Furthermore, we can approximate locally $\bm{h}_s(\bm{r},t)$ up to degree $p$ in space, by choosing orthogonal polynomials $\hat{\boldsymbol{\varphi}}_{j\ell}$ (with $j \leq d-1$) on a reference facet $\hat{F}=\hat{T}_{d-1}$ (i.e. a segment for $d=2$) and using the same transformation rules employed for the $\boldsymbol{\varphi}_{j\ell}^T$ to obtain appropriate basis functions ${\boldsymbol{\varphi}}_{j\ell}^F$ in the global coordinates.
The differential equation for the jump variable in (\ref{eq:jmp_var_def}) will thus also be required to hold in the weak sense:
\begin{align}
& \sum\limits_{F\in\mathcal{T}_\Omega} {\alpha}^{-1} \left(\mu \partial_t \bm{h}_s,\boldsymbol{\varphi}_{j\ell}^F\right)_{F} = \sum\limits_{F\in\mathcal{T}_\Omega} \left([[\bm{e}]]_{\hat{\bm{t}}_F}, \boldsymbol{\varphi}_{j\ell}^F\right)_{F} ,  \label{eq:weak_penalty}
\end{align} for all ${\boldsymbol{\varphi}}_{j\ell}^F\in\mathbb{P}(F;\mathbb{R}^{d-1})$.
The final structure of the system is only slightly changed (although more unknowns are present) with respect to the usual one for ME, as the semi-discrete formulation for (\ref{eq:cwc}), (\ref{eq:1st_weak_equation}) and (\ref{eq:weak_penalty}) admits the following matrix representation
 \begin{align*}
              \begin{pmatrix*}[c]
               \mathbf{M}_{\varepsilon} & \mathbf{0} & \mathbf{0} \\
               \mathbf{0} &  \mathbf{M}_\mu &  \mathbf{0} \\
               \mathbf{0} & \mathbf{0} & \mathbf{M}_{\mu/\alpha} 
              \end{pmatrix*}\frac{\mathrm{d}}{\mathrm{d}t}
              \begin{pmatrix*}[c]
               {\mathbf{e}} \\
               {\mathbf{h}} \\
               \mathbf{h}_s
              \end{pmatrix*}\!=\!
              \begin{pmatrix*}[c]
               \mathbf{0} & \mathbf{C} & \mathbf{C}_{F} \\
               -\mathbf{C}_{}^{\mathrm{T}} & \mathbf{0} & \mathbf{0} \\
               -\mathbf{C}_{F}^{\mathrm{T}} & \mathbf{0} & \mathbf{0}
              \end{pmatrix*}
              \begin{pmatrix*}[c]
               {\mathbf{e}} \\
               {\mathbf{h}} \\
               \mathbf{h}_s
              \end{pmatrix*},
 \end{align*}
\noindent where we have grouped DoFs in column vectors: $\mathbf{e}$, $\mathbf{h}$ are the column vectors containing time-dependent DoFs $u_{j\ell}(t)$, $f_{j\ell}(t)$ respectively, while $\mathbf{h}_s$ has analogous role for the penalty field. We refer the reader to \cite{warburton_role_2006} and \cite{koutschan_computer_2012} for more thorough discussions on the properties of the penalization term.
\section{Perfectly Matched layers}\label{sec:pml_tf}
Scattering problems present solutions for which the energy density decays very slowly as they radiate towards $\bm{r}\!\rightarrow\! +\infty$, yet any computer can only store the mesh of a bounded domain. It is therefore a key feature for our numerical method to succeed in truncating unbounded domains with perfectly absorbing (i.e. also non-reflecting) boundary conditions. Since the work of Berenger\cite{berenger_perfectly_1994}, the Perfectly Matched layer (PML) has become a standard choice in the literature. A fruitful formulation of the PML as a complex-valued stretching of the Cartesian spatial coordinates in $\mathbb{R}^d\setminus\Omega$ was introduced in \cite{chew_3d_1994}. For ease of presentation, we will from now on assume $\Omega$ to be a cuboid:
\begin{align*}
&\Omega= \big\{ \bm{r}(x,y,z) : |x|\!\leq\!x^*, |y|\!\leq\!y^*, |z|\!\leq\!z^*\big\},
\end{align*}
\noindent with $x^*$, $y^*$, $z^*$ positive real numbers.
The basic idea consists in enclosing the original computational domain with its complement $\Omega^*$ with respect to a bigger cuboid $\tilde{\Omega}=\Omega\cup\Omega^*$. For conformal PMLs in more general geometries we refer the reader to \cite{teixeira_systematic_1997,collino_optimizing_1998}.
In $\Omega^*$ we can employ a PML which consists of separately stretching each global Cartesian coordinate of space by a complex-valued factor, i.e.
\begin{align}
  & \tilde{\bm{r}} = \tilde{\Phi}_T(\bm{r}) = \mathbf{I}\bm{r} + i
  \begin{pmatrix*}[c] \alpha_x & 0 & 0\\
    0 & \alpha_y & 0\\
    0 & 0 & \alpha_z \end{pmatrix*}\left(\bm{r}-\bm{r}^*\right),\label{eq:cartesian_stretching}
\end{align}
\noindent where $\mathbf{I}$ is the ${d\!\times\!d}$ identity matrix, $i=\sqrt{-1}$, $\alpha_x$,$\alpha_y$,$\alpha_z\!\in\mathbb{R}^+$, and 
\begin{align*}\bm{r}^*\!=\!\begin{pmatrix*}[c] x^* sgn(x) \\y^* sgn(y)\\ z^* sgn(z)\end{pmatrix*},\end{align*} with $sgn(\cdot)$ being the sign function. It can be easily shown, by expanding the continuous problem's solution into wave-like solutions in a homogeneous medium (as in \cite{johnson_notes_2010}), that the imaginary part of the stretching transforms oscillating solutions into exponentially decaying ones, without changing the wave-impedance of incoming waves in the process.
\begin{figure}[!t]
\centering
\includegraphics[width=\columnwidth]{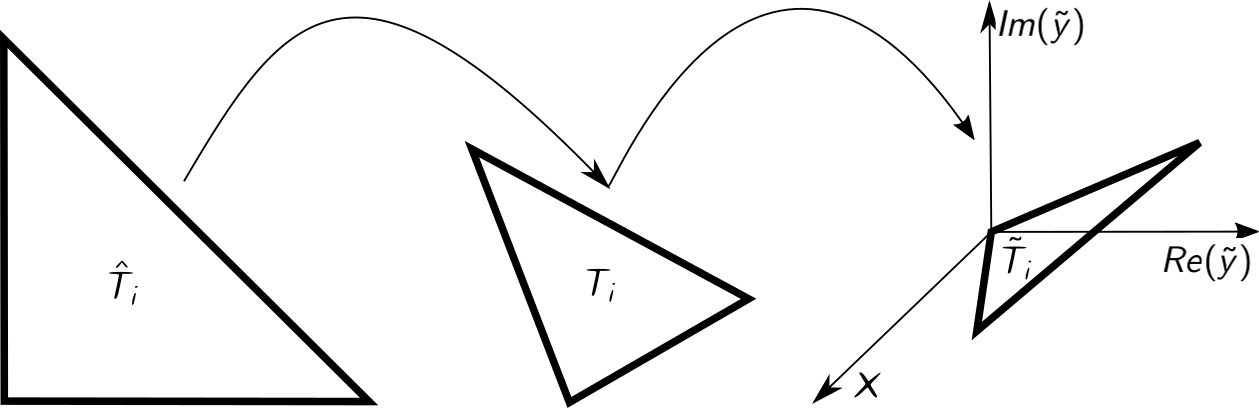}
\caption{The reference triangle, on which trial and test functions are defined in a FEM mesh, undergoes a composition of transformations in $\mathbb{C}^2$. }
\label{fig:element_maps}
\end{figure}

If we now re-introduce the triangulation $\mathcal{T}_\Omega$, one can readily notice that Eq. (\ref{eq:cartesian_stretching}) describes a family of mappings $\tilde{\bm{r}} = \tilde{\Phi}_{T}(\bm{r})$ similar to (\ref{eq:mapping_template}), where one just allows a complex valued $\mathbf{A}_{\tilde{T}}$ in lieu of $\mathbf{A}_T$. Moreover, the two element mappings ${\Phi}_{T}$ and $\tilde{\Phi}_{T}$ are used in composition, whose sketch is shown in Fig. \ref{fig:element_maps} for a single triangle.
Consequently, again using the chain-rule and the mapping $\tilde{\Phi}_{T}$, electric and magnetic fields undergo new covariant transformations from any physical element $T$, to the \emph{absorbing} (complex--stretched) one $\tilde{T} \subset \mathbb{C}^d$. 

The presence of the PML renders the (covariantly) transformed electromagnetic fields, which we label $\tilde{\bm{e}}(\tilde{\bm{r}},t)$, $\tilde{\bm{h}}(\tilde{\bm{r}},t)$, complex-valued. Naturally, complex-valued vector fields have no meaning in the time domain formulation, therefore we briefly recast the ME in the (angular) frequency-domain, where we will denote with capital calligraphic letters the Fourier transforms of vector fields, e.g.
\begin{align*}
  &\bm{E}(\bm{r},\omega) = \mathcal{F}\left(\bm{e}(\bm{r},t)\right) := \int_{-\infty}^{+\infty} \bm{e}(\bm{r},t)e^{-i\omega t} \mathrm{d}t.%\\
  %% &\bm{e}(\bm{r},t) = \mathcal{F}^{-1}\left(\bm{e}(\bm{r},t)\right) := \int_{-\infty}^{+\infty}\bm{E}(\bm{r},\omega) e^{-i\omega t} \mathrm{d}\omega.
\end{align*}

The time derivatives in the weak formulation are thus Fourier-transformed into $-i\omega$ factors. To have frequency independent absorption we let $\alpha_{x,y,z}$ depend on $\omega$: we take a single real parameter $\omega^*>0$ (with physical units of angular frequency) and we set
\begin{align*}
  \alpha_x = \begin{cases}
    0 & |x| < x^*,\\
    \omega^*/{\omega} & \text{otherwise},\end{cases}
\end{align*}
\noindent with analogous definitions valid for $\alpha_y$ and $\alpha_z$.
For the l.h.s. of the weak formulation, in the Fourier domain, it ensues
\begin{align*}
 -i \omega \left(\varepsilon(\bm{r}) \tilde{\bm{E}} , \tilde{\boldsymbol{\varphi}}_{j\ell}^T\right)_{\tilde{T}} 
&=- i \omega \left( |\mathbf{A}_{\tilde{T}}|
  \left({\mathbf{A}}_{\tilde{T}}^{-1} \varepsilon {\mathbf{A}}_{\tilde{T}}^{-\mathrm{T}}  \right) {\bm{E}} , \boldsymbol{\varphi}_{j\ell}^T\right)_T\!=\\
&= -i\omega\left(\, \tilde{\varepsilon}(\bm{r},\omega){\bm{E}} , \boldsymbol{\varphi}_{j\ell}^T\right)_T,
\end{align*}
\noindent for all $\boldsymbol{\varphi}_{j\ell}^T$, where $\tilde{\varepsilon}$ will be an anisotropic material tensor which always admits an additive decomposition of the kind:
\begin{align*}
&-i\omega\tilde{\varepsilon} = -i\omega \varepsilon + \eta\varepsilon  + \left(\frac{1}{\omega^* - i\omega}\right) \xi\varepsilon,
\end{align*}
\noindent with $\eta$, $\xi$ independent from $\omega$, but dependent on position through $\omega^*$.
Inverse-Fourier transforming the single terms yields:
\begin{align*}
-i\omega \varepsilon \bm{E}(\bm{r},\omega) &\xmapsto[]{\mathcal{F}^{-1}} \varepsilon \partial_t \bm{e}(\bm{r},t), \\
    \eta \varepsilon \bm{E}(\bm{r},\omega)   &\xmapsto[]{\mathcal{F}^{-1}} \eta \varepsilon\bm{e}(\bm{r},t),
\end{align*}
\noindent for the first two terms, while for the third term we define
\begin{align*}
&\left(\frac{1}{\omega^* - i\omega}\right) \xi\varepsilon \bm{E}(\bm{r},\omega) := \bm{P}(\bm{r},\omega),\end{align*}
\noindent which implies
\begin{align}\partial_t \bm{p}(\bm{r},t) + \omega^* \bm{p}(\bm{r},t)  = \xi\varepsilon \bm{e}(\bm{r},t).\label{eq:p_variable}
\end{align}

We note that, equivalently to what has been derived for PMLs applied to the Yee algorithm, for a general Cartesian stretching in a corner region, it holds:
\begin{align*}
    \tilde{\varepsilon}(\omega) &=
    |\mathbf{A}_{\tilde{T}}| \, \left(\mathbf{A}_{\tilde{T}}^{-1} \varepsilon \mathbf{A}_{\tilde{T}}^{-\mathrm{T}} \right) = \nonumber \\
    &=\!\begin{pmatrix*}[c]
              \frac{ \left(1 + i\alpha_y\right) \left(1 + i\alpha_z\right)}{ \left(1 + i\alpha_x\right)} \!&\!  0 \!&\! 0 \\
              0  \!&\!  \frac{ \left(1 + i\alpha_z\right) \left(1 + i\alpha_x\right)}{ \left(1 + i\alpha_y\right)} \!&\!  0 \\
              0 \!&\! 0 \!&\! \frac{ \left(1 + i\alpha_x\right) \left(1 + i\alpha_y\right)}{ \left(1 + i\alpha_z\right)} 
   \end{pmatrix*} \varepsilon,
\end{align*}

\noindent while we remark, on the other hand, that bilinear forms involving space derivatives are again independent from the particular mapping (and thus computationally inexpensive). This is a pivotal result of the present article, and also gives a more rigorous variational justification for the anisotropic absorber interpretation of the PML found in the literature: all changes in the weak formulation are absorbed in the l.h.s., specifically in the material tensor behaviour with respect to $\omega$.

Analogous considerations for the magnetic field lead to the introduction of an additional unknown $\bm{q}(\bm{r},t)$ such that:
\begin{align}
\partial_t \bm{q}(\bm{r},t) + \omega^* \bm{q}(\bm{r},t)  = \xi\mu\bm{h}(\bm{r},t),\label{eq:q_variable}
\end{align}
where additional unknowns $\bm{p}$ and $\bm{q}$ introduced in the time domain formulation can be discretized within the same function spaces as the original unknowns, therefore all considerations on efficiency due to the underlying orthogonal basis seamlessly extend to the additional DoFs relative to new unknowns. The modified system, in the new computational domain $\tilde{\Omega}$ and including weak versions of (\ref{eq:p_variable})--(\ref{eq:q_variable}), can be then solved with the usual perfect electric conductor (PEC) boundary conditions ($\bm{e}|_{\partial\tilde{\Omega}}\! \times \!\hat{\bm{n}}(\partial\tilde{\Omega}) = 0$), since the fields are exponentially damped in the PML medium.

\begin{figure}[!t]
    \centering
    \includegraphics[width=0.59\linewidth]{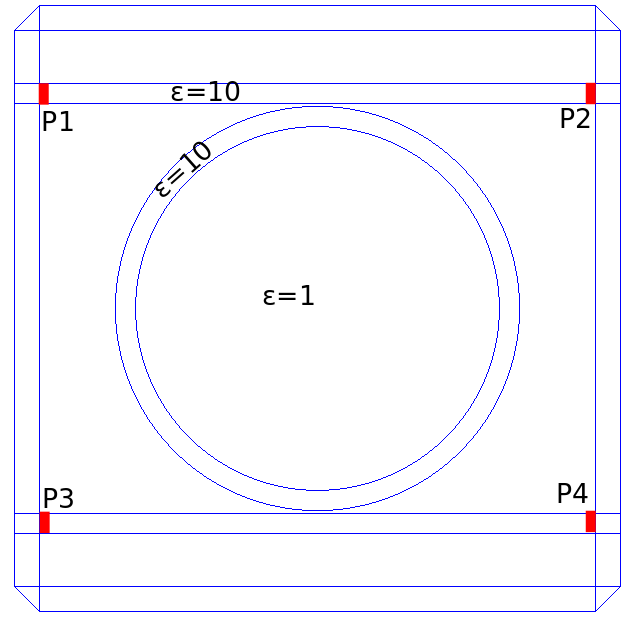}
    \caption{The test waveguide ring geometry: the PML wrap-around is also shown.}
    \label{fig:geometry}
\end{figure}

\section{Numerical Results}\label{sec:num_res}
%The main application advocated by engineers for time domain solutions of ME is the extraction of broadband frequency power spectra through a single simulation.
We shall hereinafter present numerical results which validate the proposed DG method. All experiments shown in the present section have been obtained within the Finite Element library NGSolve \cite{schoberl_c11_nodate} (built on top of the open-source 3D mesher NetGen\cite{schoberl_netgen_1997}). The library is developed in the C++ language, but has a Python front-end in which finite element spaces, algebraic system solvers and post-processing techniques can be invoked in a user-friendlier environment. We show in Appendix \ref{sec:ngcode} the Python script needed for a full three-dimensional problem employing the PML parameters and structures defined in Section \ref{sec:pml_tf}.

All tests in the present section have been computed using normalized units. We have discretized the time variable with the classic leap-frog scheme and kept the time-stepping explicit by treating all dissipative terms with the trapezoid integration rule, which is also second-order accurate in time.
We show first the results for a two-dimensional problem, namely a dielectric slab-waveguide-ring structure with four ports (P1 through P4 in Fig.~\ref{fig:geometry}). This is a common optical circuit, in which a wave-packet is injected at the port $P1$ into the upper slab waveguide, where we set (relative) permittivity $\varepsilon=10$. Here the radiation is guided by virtue of total internal reflection, until it interacts, through evanescent waves existing in the small gap, with the ring structure (having same $\varepsilon$ as the slab). The fraction of energy which is thus guided inside the ring will again couple with the lower waveguide, and a fraction of it will be guided trough port P3. We use polynomial orders $p=5$ for the $\bm{h}$ field and $p+1=6$ for the $\bm{e}$ field. Geometrical parameters for the structure are taken from \cite{busch_discontinuous_2011}. The domain is augmented by a layer of Cartesian PML with width one tenth the one of the actual $\Omega$.
We remark that in two dimensions one of the two fields is a pseudo-vector aligned with the $z$--axis. This poses no real complication and all appropriate modifications in the derivations of the previous sections are straightforward.

\begin{figure}[!t]
    \centering
    \includegraphics[width=0.8\linewidth]{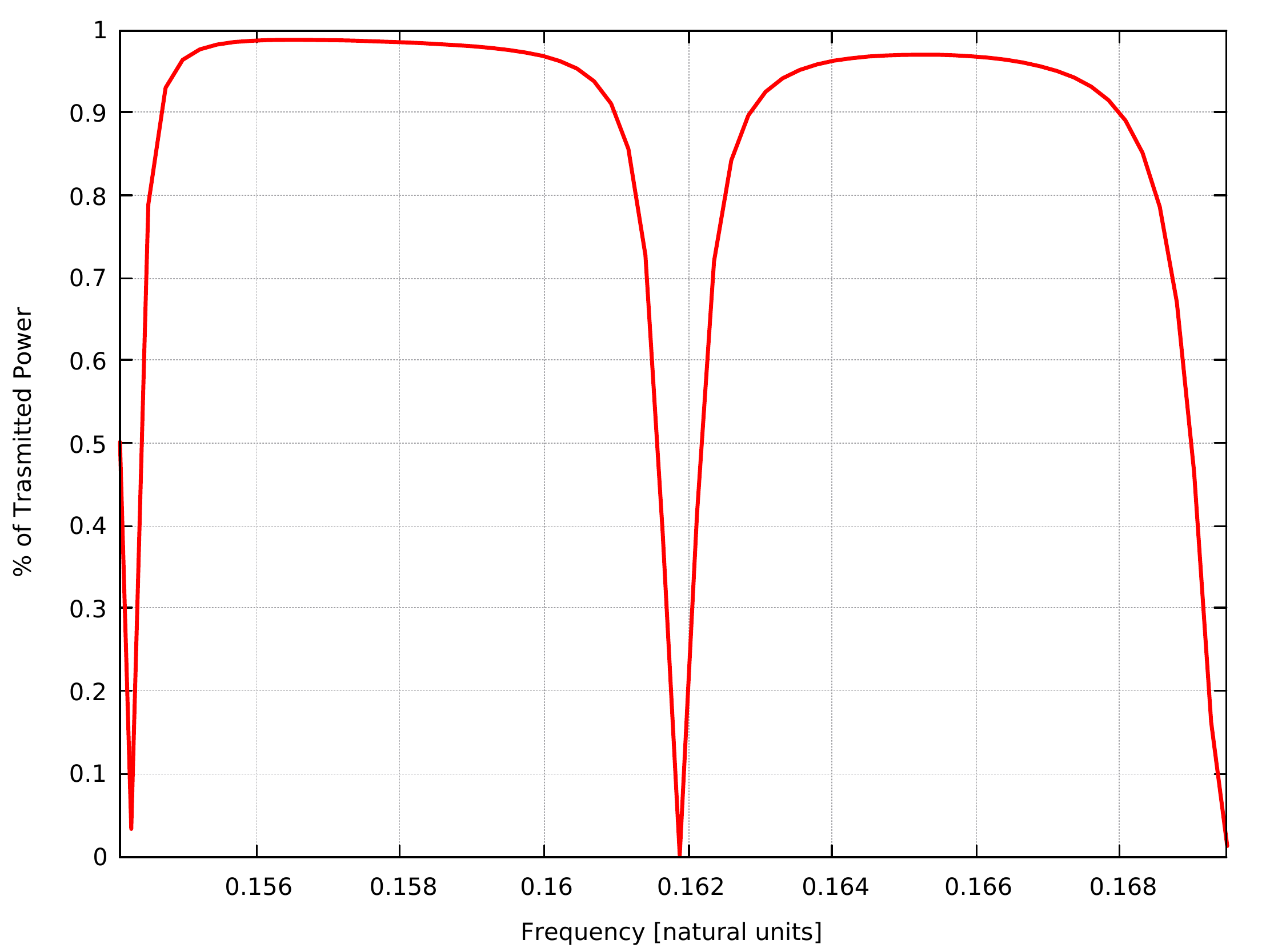}
    \caption{Ratio between power absorbed from P2 over power injected at P1, obtained by FFT.}
    \label{fig:s12}
\end{figure}

\begin{figure}[!t]
    \centering
    \includegraphics[width=0.6\linewidth]{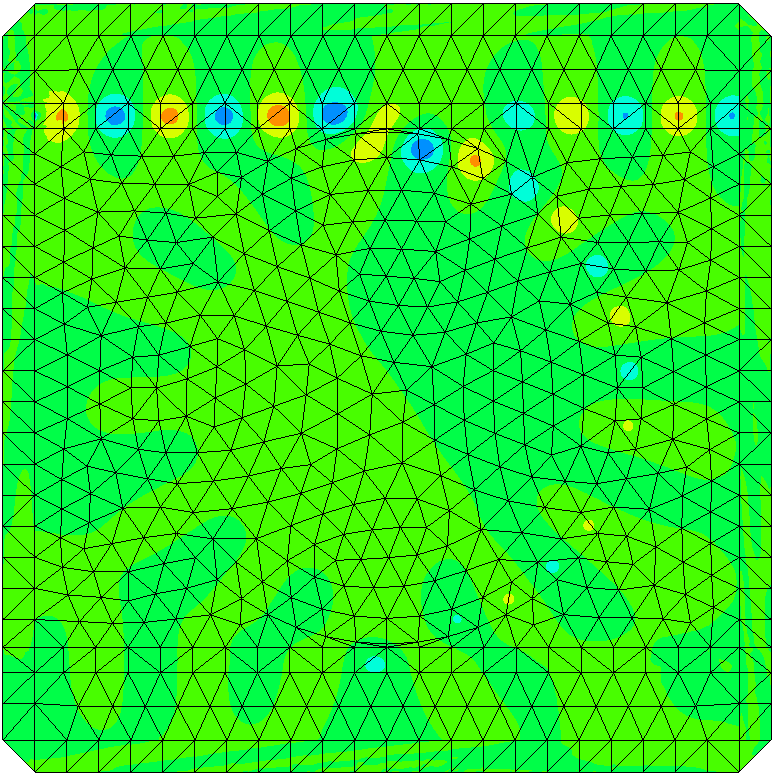}
    \caption{First few wave-fronts of the time domain simulation.}
    \label{fig:transient}
\end{figure}

\begin{figure}[!t]
    \centering
    \includegraphics[width=0.6\linewidth]{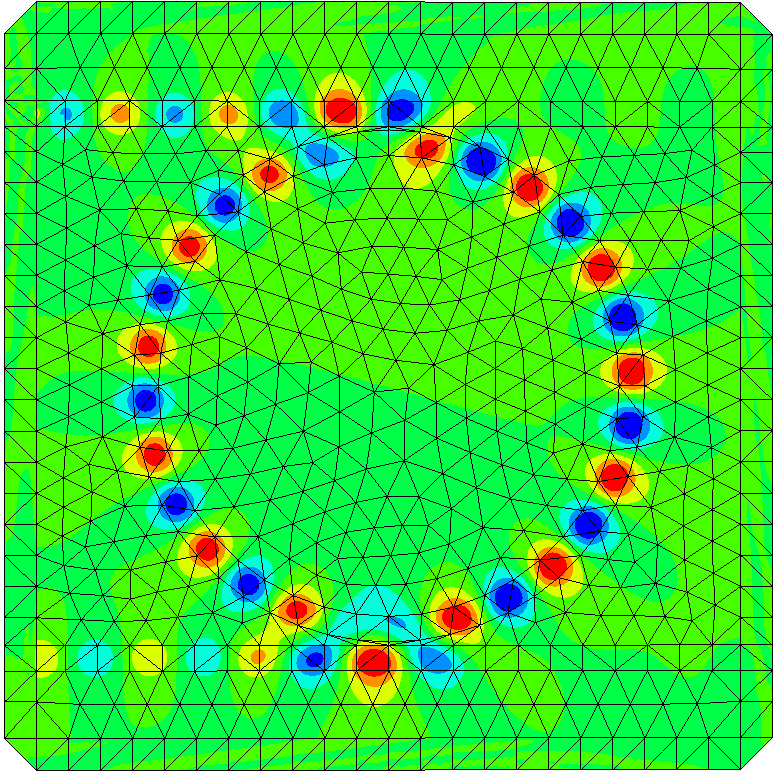}
    \caption{Steady-state fields after 400\,000 time steps.}
    \label{fig:steadystate}
\end{figure}

We begin by performing a time domain simulation in which a source with wide bandwidth (in the frequency domain) is employed. Through fast-Fourier-transform (FFT) analysis on the computed fields, the spectrum of our test structure can be studied in terms of transmission coefficient, e.g.
\begin{align}
\text{S}_{\text{21}} := \frac{|\int_{\text{P2}} \bm{E}(\bm{r},\omega)\!\times\!\bm{H}(\bm{r},\omega)\cdot\hat{\bm{n}}(\text{P2})\mathrm{d}S|}
                             {|\int_{\text{P1}} \bm{E}(\bm{r},\omega)\!\times\!\bm{H}(\bm{r},\omega)\cdot\hat{\bm{n}}(\text{P1})\mathrm{d}S|},\label{eq:s21}
\end{align}
\noindent for which a snapshot with respect to $\omega$ is shown in Fig.~\ref{fig:s12}, which reveals the expected filtering behaviour. The integrals in (\ref{eq:s21}) can be computed via Gauss quadrature rules on FFT post-processed samples of the fields on appropriate segments of the mesh.
We can validate the result of Fig.~\ref{fig:s12} a posteriori, by injecting a pure sine mode whose frequency matches a computed resonant one: as shown in Fig.\ref{fig:transient}, the energy of the sinusoidal field injected in the slab waveguide is initially split, through evanescent mode coupling, between the ring structure and the waveguide's upper-right output port.
After a full round-trip around the ring structure, the guided mode interacts (again through evanescent mode coupling) destructively with the continuous wave coming from the source.
The ensuing steady state behaviour comprises resonant energy stored in the ring structure, and a guided mode from the P1 to P3, shown in Fig.~\ref{fig:steadystate}.

\begin{figure}[!t]
    \centering
    \includegraphics[width=0.8\columnwidth]{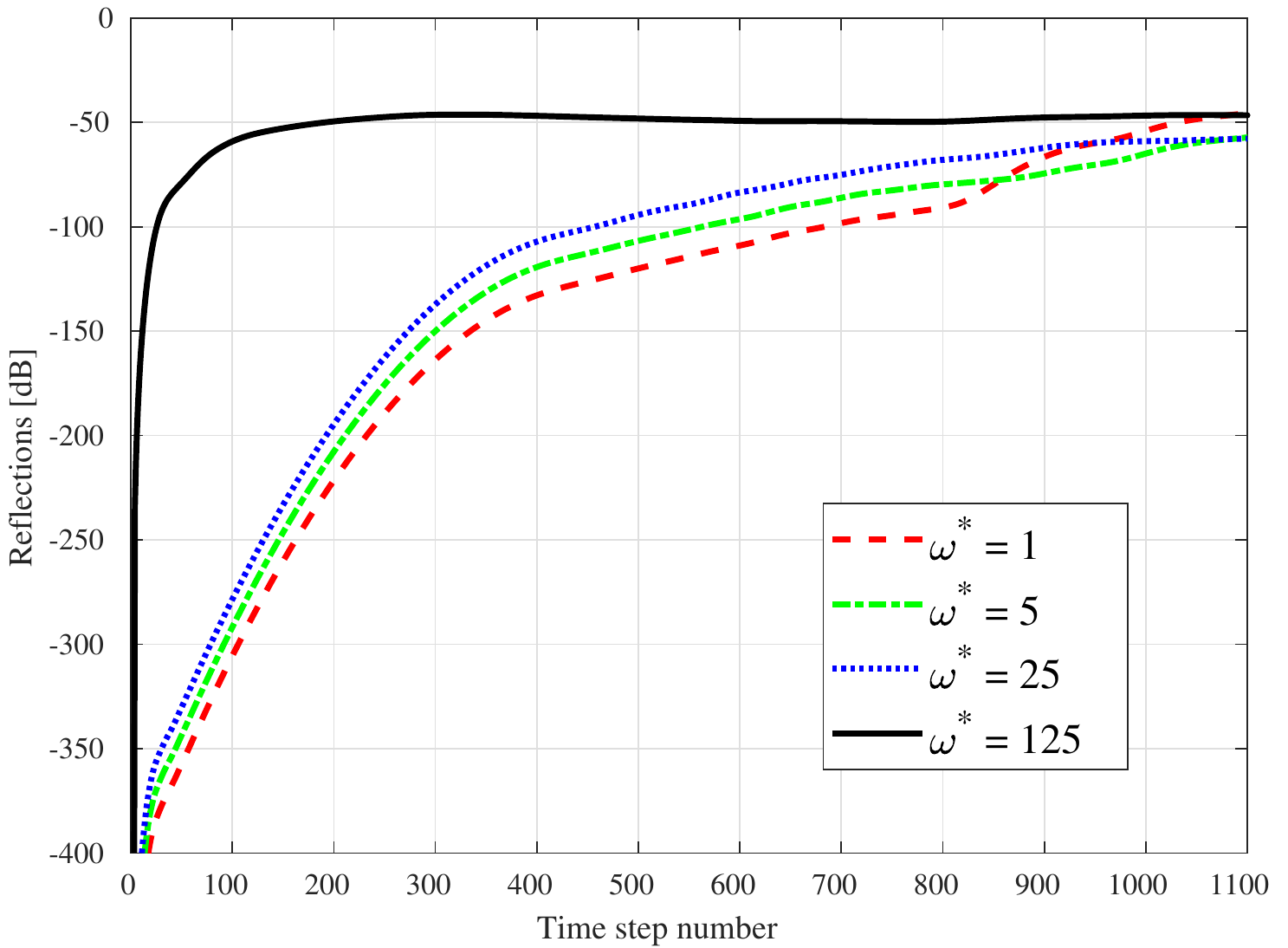}
    \caption{ The $\mathbf{L}^2$ norm of the reflections for the scattered field formulation example. 
        The error is computed with respect to another simulation with the same mesh and excitation, in which PMLs anisotropic $\omega^*$ has been set to zero.}
    \label{fig:various_sigma}
\end{figure}

To test the method in three space dimensions, we compute the scattered field from a PEC sphere. The incident field is a plane wave propagating in the $z$ direction with known amplitude, modulated in time by a smooth function with compact support.
The quality of the PML can be tested in a standard way by performing a second simulation where the interior domain $\Omega$ is augmented in size by a factor two in all directions and no PML is used. Fig.~\ref{fig:various_sigma} shows the reflections caused by the discretized PML in the $L^\infty$ norm, for various values of the scaling parameter $\omega^*$, with respect to time.
The sensitivity is certainly non-negligible: low values of $\omega^*$ imply lower reflections from the interface with the PML, which are the earliest ones arising. On the other hand, low damping implies higher reflections arising from the truncation of the new domain $\tilde{\Omega}$ (we note the cusp in the red, dashed curve in Fig.~\ref{fig:various_sigma} after roughly 850 time steps).

We remark that a piecewise-uniform damping parameter $\omega^*$ has been chosen in the formulation. This is motivated by performance arguments, since all local Jacobians arising from the mappings in (\ref{eq:cov_transf}) and (\ref{eq:cartesian_stretching}) have in this case constant coefficients. No numerical integration must be then performed to obtain the discrete operators needed in the final algorithm. 
A smoother profile for $\omega^*$ would obviously lead to less reflection at the interface with the interior domain. On the other hand, accurate discretization of the PML domain would then become a much more critical issue.

\section{Conclusion}\label{sc:cr}
We have presented a DG method for the time dependent ME combining an orthogonal modal basis with the appropriate transformation rules, under changes of coordinates, for the physical fields involved. This combination is fruitful since it yields a very low storage implementation (competitive with the one produced for the nodal basis in \cite{warburton_low-storage_2013} or for the Cell Method in \cite{kapidani_time-domain_2020,kapidani_arbitrary-order_2020}) which additionally has only moderate increase in the number of unknowns when a PML is used to truncate an unbounded computational domain.
The optimization of our PML and its efficient generalization to arbitrary damping profile and curved elements is a current object of study and will be reported elsewhere.

\section*{Acknowledgements}
Author Bernard Kapidani is supported by the Austrian Fund for Scientific Research (FWF) under grant number F65, project: ``Taming Complexity in Partial Differential Equations''.

\begin{figure}[!t]
    \centering
    \includegraphics[width=0.8\columnwidth]{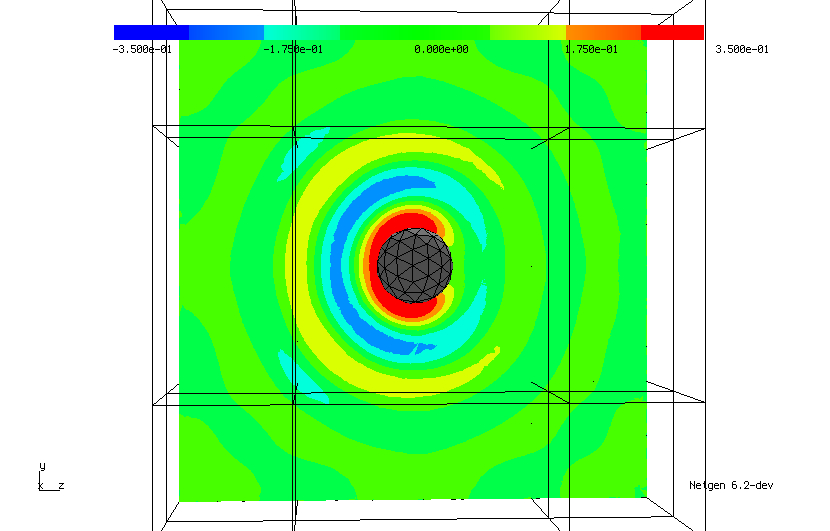}
    \caption{ A snapshot of the reflected field from the sphere on a cutting plane of the 3D test domain: the field vanishes with negligible reflections once it enters the PML.}
    \label{fig:snapshot}
\end{figure}

\begin{appendices}
\section{NGSolve listing}\label{sec:ngcode}

In the following we present the script needed to simulate the full 3D problem discussed in Section \ref{sec:num_res}, inclusive of PML. For ease of presentation we leave out the interior penalty terms.
After importing the appropriate libraries, NetGen allows the definition of geometry and meshing parameters of the domain $\Omega$ through Constructive Solid Geometry (CSG) instructions, starting from the background cuboid with unit length sides:
{\tiny\begin{Verbatim}[frame=none]
#################################################################################
from ngsolve import *
from netgen.csg import *
from netgen import gui
SetHeapSize = 10*100*1000
# half-side length of the domain
xi = 0.5
yi = 0.5
zi = 0.5
xo = 0.6
yo = 0.6
zo = 0.6
geo  = CSGeometry()
omega  = (OrthoBrick(Pnt(-xi,-yi,-zi),Pnt(xi,yi,zi))) #original domain
#################################################################################
\end{Verbatim}
}
\noindent All the other cuboids (``\lstinline|Orthobrick|'') pertaining to the PML can be defined in similar fashion. Finally, only the domains which actually need to be meshed are added to the ``\lstinline|geo|'' object and the graphical user interface (GUI) is asked to show the mesh via a ``\lstinline|Draw|'' command.
{\tiny\begin{Verbatim}[frame=none]
#################################################################################
xpml  = (OrthoBrick(Pnt(-xo,-yi,-zi),Pnt(xo,yi,zi))-omega).mat("xpml")
ypml  = (OrthoBrick(Pnt(-xi,-yo,-zi),Pnt(xi,yo,zi))-omega).mat("ypml")
zpml  = (OrthoBrick(Pnt(-xi,-yi,-zo),Pnt(xi,yi,zo))-omega).mat("zpml")
xypml = (OrthoBrick(Pnt(-xo,-yo,-zi),Pnt(xo,yo,zi))-omega-xpml-ypml).mat("xypml")
xzpml = (OrthoBrick(Pnt(-xo,-yi,-zo),Pnt(xo,yi,zo))-omega-xpml-zpml).mat("xzpml")
yzpml = (OrthoBrick(Pnt(-xi,-yo,-zo),Pnt(xi,yo,zo))-omega-ypml-zpml).mat("yzpml")
all_other_pml = xpml+ypml+zpml+xypml+xzpml+yzpml+omega
xyzpml = (OrthoBrick(Pnt(-xo,-yo,-zo),Pnt(xo,yo,zo))-all_other_pml).mat("xyzpml")
hole   = (Sphere(Pnt(0,0,0),0.1)).bc("incidentfield")
omega  = (omega-hole).mat("air") #new air filled domain

geo.Add(omega)
geo.Add(xpml)
geo.Add(ypml)
geo.Add(zpml)
geo.Add(xypml)
geo.Add(xzpml)
geo.Add(yzpml)
geo.Add(xyzpml)
mesh = Mesh(geo.GenerateMesh(maxh=0.1)) # maxh sets max mesh-size
Draw(mesh)
#################################################################################
\end{Verbatim}
}
\noindent The function ``\!\lstinline|CoefficientFunction|\!'' allows the piecewise-smooth definition of scalar, vector and matrix valued functions on $\Omega$ or $\partial\Omega$.
This is useful, for example, in defining the standard material tensors, which are $3\times3$ matrices. We define also additional useful parameters in the simulation and the special function ``\lstinline|n|'', which is the normal unit vector on each finite element boundary, and it is always outwards pointing.
{\tiny\begin{Verbatim}[frame=none]
#################################################################################
epsr  = CoefficientFunction( (1, 0, 0, 0, 1, 0, 0, 0, 1), dims=(3,3))
mur   = CoefficientFunction( (1, 0, 0, 0, 1, 0, 0, 0, 1), dims=(3,3))
mu0 = 1
eps0 = 1
n = specialcf.normal(mesh.dim)
#################################################################################
\end{Verbatim}
}
\noindent The additional fictitious material tensors $\eta(\bm{r})$, $\xi(\bm{r})$ required for the PML are also defined piecewise through dictionaries:
{\tiny\begin{Verbatim}[frame=none]
#################################################################################
sig = 5 # sig=0 => no pml
eta = {}
xi = {}

eta["air"]    = CoefficientFunction((0,0,0,0,0,0,0,0,0),dims=(3,3))
eta["xpml"]   = CoefficientFunction((-sig,0,0,0,sig,0,0,0,sig),dims=(3,3))
eta["ypml"]   = CoefficientFunction((sig,0,0,0,-sig,0,0,0,sig),dims=(3,3))
eta["zpml"]   = CoefficientFunction((sig,0,0,0,sig,0,0,0,-sig),dims=(3,3))
eta["xypml"]  = CoefficientFunction((0,0,0,0,0,0,0,0,2*sig),dims=(3,3))
eta["xzpml"]  = CoefficientFunction((0,0,0,0,2*sig,0,0,0,0),dims=(3,3))
eta["yzpml"]  = CoefficientFunction((2*sig,0,0,0,0,0,0,0,0),dims=(3,3))
eta["xyzpml"] = CoefficientFunction((sig,0,0,0,sig,0,0,0,sig),dims=(3,3))

xi["air"]    = CoefficientFunction((0,0,0,0,0,0,0,0,0),dims=(3,3))
xi["xpml"]   = CoefficientFunction((sig**2,0,0,0,0,0,0,0,0 ),dims=(3,3))
xi["ypml"]   = CoefficientFunction((0,0,0,0,sig**2,0,0,0,0 ),dims=(3,3))
xi["zpml"]   = CoefficientFunction((0,0,0,0,0,0,0,0,sig**2),dims=(3,3))
xi["xypml"]  = CoefficientFunction((0,0,0,0,0,0,0,0,sig**2),dims=(3,3))
xi["xzpml"]  = CoefficientFunction((0,0,0,0,sig**2,0,0,0,0),dims=(3,3))
xi["yzpml"]  = CoefficientFunction((sig**2,0,0,0,0,0,0,0,0),dims=(3,3))
xi["xyzpml"] = CoefficientFunction((0,0,0,0,0,0,0,0,0),dims=(3,3))

eta = CoefficientFunction([eta[mat] for mat in mesh.GetMaterials()])
xi  = CoefficientFunction([xi[mat]  for mat in mesh.GetMaterials()])
#################################################################################
\end{Verbatim}
}
\noindent The finite element spaces are then defined, for fields $e$, $h$ and auxiliary unknowns $p$, $q$. 
The keyword argument ``\lstinline|covariant|'' is fundamental, since it specifies that, even if the spaces are fully discontinuous (``\lstinline|VectorL2|''), they are mapped from reference to physical element with the appropriate transformation.
Keyword ``\lstinline|all_dofs_together|'' forces DoFs pertaining to each element to be adjacent in the global vector.
Unknowns are then grouped together into a Cartesian product space for convenience, while all their the DoFs (stored in ``\lstinline|solution|'') are initially set to zero and plotted in the GUI:
{\tiny\begin{Verbatim}[frame=none]
#################################################################################
degree = 6
fes_e = VectorL2(mesh, order=degree+1,covariant=True, all_dofs_together=True)
fes_h = VectorL2(mesh, order=degree,  covariant=True, all_dofs_together=True)
fes_p = VectorL2(mesh, order=degree+1,covariant=True, all_dofs_together=True)
fes_q = VectorL2(mesh, order=degree,  covariant=True, all_dofs_together=True)
fes   = FESpace( [fes_e,fes_h,fes_p,fes_q] )

edofs = fes.Range(0)
hdofs = fes.Range(1)
pdofs = fes.Range(2)
qdofs = fes.Range(3)

sol = GridFunction(fes)
sol.components[0].vec[:] = 0
sol.components[1].vec[:] = 0
sol.components[2].vec[:] = 0
sol.components[3].vec[:] = 0
#################################################################################
\end{Verbatim}
}
\noindent To let operators act on the global vector of DoFs and on the spaces for single unknowns without creating too many copy of involved objects, we define embedding operators, which amount to diagonal, positive-semidefinite matrices (containing only zeros and ones), acting on the DoFs vector. Operator composition is achieved through the ``\lstinline|@|'' token.
{\tiny\begin{Verbatim}[frame=none]
#################################################################################
emb_e = Embedding(fes.ndof, fes.Range(0))
emb_h = Embedding(fes.ndof, fes.Range(1))
emb_p = Embedding(fes.ndof, fes.Range(2))
emb_q = Embedding(fes.ndof, fes.Range(3))
#################################################################################
\end{Verbatim}
}
\noindent Thus, after computing the all needed mass matrices (including necessary damping terms), we can embed them in operators acting on the global solution vector:
{\tiny\begin{Verbatim}[frame=none]
################################################################################
tau    = 0.001 #time step size
alpha  = 0.5*sig*tau
beta   = 1/(1+alpha)
gamma  = (1-alpha)*beta
tpar   = Parameter(0)

massmate_eta = fes_e.Mass(eps0*(epsr-0.5*tau*eta*epsr))
massmath_eta = fes_h.Mass(mu0*(mur-0.5*tau*eta*mur))
massmate_xi  = fes_e.Mass(eps0*(eta*epsr))
massmath_xi  = fes_h.Mass(mu0*(eta*mur))
invmasse_eta = fes_e.Mass(eps0*(epsr+0.5*tau*eta*epsr)).Inverse()
invmassh_eta = fes_h.Mass(mu0*(mur+0.5*tau*eta*mur)).Inverse()

mate_eta = emb_e @ massmate_eta @ emb_e.T
math_eta = emb_h @ massmath_eta @ emb_h.T
inve_eta = emb_e @ invmasse_eta @ emb_e.T
invh_eta = emb_h @ invmassh_eta @ emb_h.T
mate_xi  = emb_p @ massmate_xi @ emb_e.T
math_xi  = emb_q @ massmath_xi @ emb_h.T
#################################################################################
\end{Verbatim}
}
\noindent A Python function defines the main loop of the script, where time-stepping is actually performed: the operators for the discrete curl are still abstract and will be defined in a later code-block.
{\tiny\begin{Verbatim}[frame=none]
#################################################################################
def Run(CE, CH, t0 = 0, tend = 1, tau = 1e-3):
  t = 0
  dummy = sol.vec.CreateVector()
  with TaskManager():
    print("t = ",t)
    tpar = t
    while t < (tend-t0) - tau/2:
      t += tau
      # E, H solution update
      dummy.vec.data = inve_eta @ mate_eta*sol.vec + tau*inve_eta @ CE*sol.vec
      sol.vec[edofs] = dummy.vec[edofs]
      dummy.vec.data = invh_eta @ math_eta*sol.vec + tau*invh_eta @ CH*sol.vec
      sol.vec[hdofs] = dummy.vec[hdofs]
      Redraw(nonblocking=True)
      # auxiliary variables update
      dummy.vec.data = gamma*emb_p @ emb_p.T*sol.vec  + beta*mate_xi*sol.vec
      sol.vec[pdofs] = dummy.vec[pdofs]
      dummy.vec.data = gamma*emb_q @ emb_q.T*sol.vec  + beta*math_xi*sol.vec
      sol.vec[qdofs] = dummy.vec[qdofs]
#################################################################################
\end{Verbatim}
}
\noindent The incident field can be then enforced, via a python dictionary, on the appropriate boundary surfaces in the mesh. Namely the one labelled above with the ``\lstinline|incidentfield|'' boundary condition. 
{\tiny\begin{Verbatim}[frame=none]
#################################################################################
Einc = {}

Einc["default"] = CoefficientFunction( (0,0,0) )
Einc["incidentfield"] = CoefficientFunction( (exp(-((tpar-300*tau)/(10*tau))**2),0,0) )
Einc = CoefficientFunction([Einc[bb] for bb in (mesh.GetBoundaries())])
#################################################################################
\end{Verbatim}
}
\noindent All the terms described in Section~\ref{sec:weak_form} for the curl operator are given below in terms of bilinear-form integrators, where the syntax is inherited from the language of differential forms and the ``\lstinline|element_boudary|'' keyword denotes integration on $\partial T$, for each element $T\in\mathcal{T}_\Omega$.
Another crucial keyword argument introduced in the following code snippet is ``\lstinline|geom_free|'': this flag, when set to true, signals to the interpreter that the associated bilinear form is independent from the geometry of any particular finite element, which prompts the library to optimize integration, by building a single local reference operator (matrix).
A caveat must be mentioned: the responsibility of checking the mathematics to ensure that this is really the case falls on the user.
{\tiny\begin{Verbatim}[frame=none]
#################################################################################
h = fes_h.TestFunction()
E = fes_e.TrialFunction()

C_el  = BilinearForm(trialspace=fes_e, testspace=fes_h, geom_free = True)
C_el += -curl(E)*h*dx
C_el.Assemble()

C_tr  = BilinearForm(trialspace=fes_e, testspace=fes_h, geom_free = True)
C_tr += 0.5*Cross(E.Other(bnd=2*Einc+E)-E,n)*h*dx(element_boundary=True)
C_tr.Assemble()

C_EH = emb_h @ (C_el.mat + Ctr.mat) @ emb_e.T
C_e  =  C_EH.T - (emb_e @ emb_p.T)
C_h  = -C_EH   + (emb_h @ emb_q.T)
#################################################################################
\end{Verbatim}
}
\noindent We are finally ready to run the simulation: we tell the GUI to plot the initial (zeroed) solutions for the electro-magnetic field and we call the ``\lstinline|Run|'' function with the appropriate arguments. The ``\lstinline|Redraw|'' calls inside the loop, will update the plots at each time step.
{\tiny\begin{Verbatim}[frame=none]
#################################################################################
Draw(sol.components[1], mesh, "H")
Draw(sol.components[0], mesh, "E")
Run(C_e, C_h,0,1.2,tau)
#################################################################################
\end{Verbatim}
}
\noindent More examples and tutorials for all the basic functions in the library can be found online at \cite{noauthor_netgenngsolve_nodate}.
\end{appendices}

\ifCLASSOPTIONcaptionsoff
  \newpage
\fi
%\newpage
\bibliographystyle{./IEEEtran}
\bibliography{./IEEEabrv,./BibDB}

%\begin{IEEEbiography}[]{Bernard~Kapidani}
%received the Ph.D. degree (Cum Laude, with focus on Computational Electromagnetics) in 2018 from the University of Udine, Italy. He is currently a post-doctoral researcher at the Institute for Analysis and Scientific Computing at the Technical University of Vienna, Austria. His research interests include time domain methods for the numerical solution of Maxwell's equations and the application of computational topology and computational geometry to the solution of problems in electromagnetism. %He was short-listed as one of four finalists for the \emph{Rita Trowbridge Prize} Best Paper Award (for papers presented by a young researcher) at COMPUMAG 2017 (Daejeon, Korea).
%\end{IEEEbiography}
%\begin{IEEEbiography}[]{Joachim~Sch{\"o}berl}
%Bio goes here.
%\end{IEEEbiography}
\end{document}